\documentclass[11pt, a4paper]{amsart}

\usepackage{t1enc}
\usepackage[latin1]{inputenc}
\usepackage[english]{babel}
\usepackage{amsmath}
\usepackage{amsthm}
\usepackage{amsfonts}
\usepackage{latexsym}
\usepackage[dvips]{graphicx}
\usepackage{float}
\usepackage[natural]{xcolor}
\usepackage{algorithm}
\usepackage{algorithmic}
\usepackage{enumerate}
\usepackage{multirow}
\usepackage{xcolor}
\usepackage[colorlinks,linkcolor=blue]{hyperref}
\usepackage{lineno}
\usepackage{setspace}
\usepackage{fullpage}
\usepackage{bm}
\usepackage{verbatim}

\usepackage{listings}
\theoremstyle{plain}
\newtheorem{thm}{Theorem}[section]
\newtheorem{theorem}[thm]{Theorem}
\newtheorem{conjecture}[thm]{Conjecture}
\newtheorem{lemma}[thm]{Lemma}
\newtheorem{corollary}[thm]{Corollary}
\newtheorem{proposition}[thm]{Proposition}

\theoremstyle{definition}

\newtheorem{question}[thm]{Question}

\newtheorem{remark}[thm]{Remark}

\newtheorem{definition}[thm]{Definition}
\newtheorem{claim}[thm]{Claim}
\newtheorem{fact}[thm]{Fact}

\newtheorem{example}[thm]{Example}

\newtheorem{defn-thm}[thm]{Definition-Theorem}

\newcommand{\btheorem}{\begin{theorem}}
\newcommand{\etheorem}{\end{theorem}}
\newcommand{\bconjecture}{\begin{conjecture}}
\newcommand{\econjecture}{\end{conjecture}}
\newcommand{\bproposition}{\begin{proposition}}
\newcommand{\eproposition}{\end{proposition}}
\newcommand{\bdefinition}{\begin{definition}}
\newcommand{\edefinition}{\end{definition}}
\newcommand{\bcorollary}{\begin{corollary}}
\newcommand{\ecorollary}{\end{corollary}}
\newcommand{\bproof}{\begin{proof}}
\newcommand{\eproof}{\end{proof}}
\newcommand{\bclaim}{\begin{claim}}
\newcommand{\eclaim}{\end{claim}}
\newcommand{\bquestion}{\begin{question}}
\newcommand{\equestion}{\end{question}}
\newcommand{\bfact}{\begin{fact}}
\newcommand{\efact}{\end{fact}}
\newcommand{\bremark}{\begin{remark}}
\newcommand{\eremark}{\end{remark}}
\newcommand{\eexample}{\end{example}}
\newcommand{\bexample}{\begin{example}}

\newcommand{\elemma}{\end{lemma}}
\newcommand{\blemma}{\begin{lemma}}

\newcommand{\beq}{\begin{equation}}
\newcommand{\eeq}{\end{equation}}

\newcommand{\mathify}[1]{\ifmmode{#1}\else\mbox{$#1$}\fi}
\newcommand{\bigO}O

\newcommand{\remove}[1]{{}}

\title{On the thresholds of degenerate hypergraphs}

\author{Yu Chen}
\address{School of Mathematics and Statistics\\ Beijing Institute of Technology\\ Beijing\\ China}
\email{chenyu2023@bit.edu.cn}

\author{Jie Han}
\address{School of Mathematics and Statistics\\ Beijing Institute of Technology\\ Beijing\\ China}
\email{han.jie@bit.edu.cn}

\author{Haoran Luo}
\address{Department of Mathematics\\ University of Illinois Urbana-Champaign\\ Illinois\\ United States}
\email{haoranl8@illinois.edu}

\begin{document}
\begin{abstract}
An $n$-vertex $k$-uniform hypergraph $G$ is $(d,\alpha)$-\emph{degenerate} if $m_1(G)\le{d}$ and there exists a constant $\varepsilon >0$ such that for every subset $U\subseteq{V(G)}$ with size $2\le|U|\le{\varepsilon n}$, we have $e\left(G[U]\right)\le{d\left(|U|-1\right)-\alpha}$.
These hypergraphs include many natural graph classes, such as the degenerate hypergraphs, the planar graphs, and the power of cycles.
In this paper, we consider the threshold of the emergence of a $(d,\alpha)$-degenerate hypergraph with bounded maximum degree in the Erd\H{o}s--R\'enyi model.
We show that its threshold is at most $n^{-1/d}$,
improving previous results of Riordan and Kelly--M{\"u}yesser--Pokrovskiy.

\end{abstract}
\maketitle

\section{Introduction} \label{sec::int}
Since its inception by Erd\H os and R\'enyi in 1960~\cite{erdos1960evolution}, the random graph model has been one of the main objects of study in Probabilistic Combinatorics.
Given positive integers $n,k$ and a real number $p\in[0,1]$, $G^{(k)}(n,p)$ is the random $k$-uniform hypergraph on $n$ vertices, where each $k$-set is included as an edge with probability $p$ independently.
A central topic in random graph theory is to determine the threshold for the property of the emergency of a graph.
More precisely, given a sequence of $k$-uniform hypergraphs $\{G_n\}$, let ${\rm{th}}(G_n)$ be a \emph{threshold function}\footnote{With this definition, constant factors in the expression of ${\rm{th}}(G_n)$ do not matter and we usually omit them.} of $(G_n)$,
that is, if $p(n)/{\rm{th}}(G_n)\to \infty$, then w.h.p.~$G^{(k)}(n,p(n))$ contains a copy of $G_n$; if $p(n)/{\rm{th}}(G_n)\to 0$, then w.h.p.~$G^{(k)}(n,p)(n)$ contains no copies of $G_n$.

In~\cite{kahn2007thresholds}, Kahn and Kalai proposed the following intriguing conjecture. For an arbitrary graph $G$, let $p_E(G)$ be the smallest $p$ such that
\[\frac{v(F)!}{|\textrm{Aut}(F)|}p^{e(F)}\ge 1\text{ for all }F\subseteq{G}.\]
Note that the quantity in the above inequality is just the expected number of copies of $F$ in $G^{(k)}(n,p)$.

\begin{conjecture}\label{conj:conj_seocnd_kk}
There is a universal constant $K$ such that
for any sequence $\{G_n\}$, a threshold function for the appearance of $\{G_n\}$ in $G(n,p)$ is at most $p_E(G_n)\log|v(G_n)|$.
\end{conjecture}

Johansson, Kahn, and Vu~\cite{Johansson2008FactorsIR} proved Conjecture~\ref{conj:conj_seocnd_kk} for $F$-factors when $F$ is a (fixed) strictly balanced graph and their result is also extended to hypergraphs.
It is natural to ask more general problems: let $\mathcal{F}(n,\Delta)$ be the collection of all $n$-vertex graphs with maximum degree at most $\Delta$, what is the threshold behavior for the appearance of an arbitrary graph $G\in\mathcal{F}(n,\Delta)$ in $G(n,p)$?
Towards the family $\mathcal{F}(n,\Delta)$, Ferber, Luh, and Nguyen~\cite{Ferber2016EmbeddingLG} considered a similar problem for the family $\mathcal{F}((1-\varepsilon)n,\Delta)$ of all almost spanning graphs with maximum degree at most $\Delta$.
In~\cite{FKNP2021}, Frankston, Kahn, Narayanan, and Park improved Ferber, Luh, and Nguyen's result, and showed that for every graph $G\in\mathcal{F}(n,\Delta)$, w.h.p. $G(n,p)$ contains a copy of $G$, provided that $p=\omega(n^{-1}\log^{1/\Delta}n)^{2/(\Delta+1)}$.

For a graph $G$ and a positive integer $i$, we let\footnote{Note that this definition of $m_2(G)$ is different from the one used for the so-called random Tur\'an theorem and random Ramsey theorem (sometimes also denoted by $m_2(G)$), which is $\max \left\{\frac{e(F)-1}{v(F)-2}: F\subseteq{G}\text{ and } v(F)> 2 \right\}$, defined with respect to the expected number of copies of $G$ containing a given edge.}
\[
m_i(G)=\max \left\{\frac{e(F)}{v(F)-i}: F\subseteq{G}\text{ and } v(F)> i \right\}.
\]
\noindent
A natural class of (hyper)graphs to consider is the class of degenerate graphs.
A graph $G$ is called $d$-degenerate if every subgraph of $G$ contains a vertex of degree at most $d$ (therefore, it is enough to consider $d\in \mathbb N$).
Using the second-moment method, Riordan~\cite{Riordan2000SpanningSO} showed a general result which implies that ${\rm{th}}(G_n)\le n^{-1/m_2(G_n)}$.
Letting $\{G_n\}$ be a sequence of $d$-degenerate graphs, Riordan's result implies that if $d\ge 3$, then ${\rm{th}}(G_n)\le{n^{-1/d}}$.
However, for $d=2$, it only gives a suboptimal result ${\rm{th}}(G_n)\le{n^{-1/3}}$, as $m_2(K_3)=3$ (and $K_3$ is $2$-degenerate).
At last, 1-degenerate graphs are forests, so their threshold is $\log n/n$.

In this paper, we study the threshold behavior {of a wider class, namely,} the $(d,\alpha)$-degenerate hypergraphs with a bounded maximum degree.
\begin{definition}
Let $k,\Delta\in{\mathbb{N}}$, and~$d,\alpha>0$.
Assume~$G$~is an $n$-vertex $k$-uniform hypergraph on $V$, and we say $G$ is a \emph{$(d,\alpha)$-degenerate} hypergraph if $m_1(G)\le{d}$ and there exists a constant $0<\varepsilon<1$ such that for any subset $U\subseteq{V}$ with size $2\le|U|\le{\varepsilon n}$, we have $e\left(G[U]\right)\le{d\left(|U|-1\right)-\alpha}$.
\end{definition}

The $(d,\alpha)$-degenerate hypergraphs include several natural graph classes.
\begin{example}

\begin{enumerate}
\item For $d\ge 2$, every $d$-degenerate graph is $(d,d/2)$-degenerate, as we have $e(G)\le d(n-d)+\binom{d}{2} = d(n-1) - \binom{d}{2}$ for every $n$-vertex $d$-degenerate graph $G$.
\item Every planar graph is $(3,3)$-degenerate as $e(G)\le 3v(G) - 6$ for every planar $G$.
\item The $r$-th power of a cycle is $(r,2-1/r)$-degenerate. More generally, the $r$-th power of a $k$-uniform tight cycle is $\left(\binom{r+k-2}{k-1},1-1/k\right)$-degenerate.
\end{enumerate}
\end{example}
\noindent
Indeed, as the constraint on the number of edges is only required to subgraphs of order $o(n)$, given a $(d,\alpha)$-degenerate $k$-uniform hypergraph, one can add some edges to pairs of vertices that are ``far-away'', and still keep the $(d,\alpha)$-degeneracy.

For the thresholds of $(d,\alpha)$-degenerate hypergraphs, \cite[Proposition 1.17]{kelly2023optimal} (see Proposition \ref{prop:Kelly_edge_spread} below) implies that ${\rm{th}}(G)\le n^{-1/m_1(G)}\log n$, and since we have $m_1(G)\le d$, it follows that ${\rm{th}}(G)\le n^{-1/d}\log n$.
We first improve on this result.

\begin{theorem}\label{thm:asthm0}
Let $k,\Delta\in\mathbb{N}$,~$d,\alpha$.
Let $G$ be an $n$-vertex $(d,\alpha)$-degenerate $k$-uniform hypergraph with $\Delta(G)\le\Delta$.
Then we have ${\rm{th}}(G)\le{n^{-1/d}}$.
\end{theorem}
\noindent
Theorem \ref{thm:asthm0} follows from a more general result (Theorem \ref{thm:asthm}) whose statement is postponed to Section~\ref{sec::spr}.

Theorem~\ref{thm:asthm0} has an interesting consequence and let us first motivate it with the following definition.
Let $d\in{\mathbb{N}}$, $0<\varepsilon<1/2$, and let $G$ be an $n$-vertex graph, and for every subset $U\subseteq V(G)$, we define its edge boundary $\partial U$ as $\partial U:=\{e\in E(G):|e\cap U|=1\}$.
We say that $G$ is $(d,\varepsilon)$-\emph{locally sparse} if for every subset $U\subseteq V(G)$ of size $d\le |U|\le{\varepsilon n}$, we have $|\partial U|\ge{d+1}$.
{Note that if $G$ is a $(d,\varepsilon)$-locally sparse $d$-regular graph, then $G$ is $(d/2,1/2)$-degenerate.
Thus, we have the following result for locally sparse regular graphs.}
\begin{corollary}\label{cor:threshold_of_locally_sparse_regular_graph}
Let $d\in\mathbb{N}$ and $\varepsilon \in (0, 1/2)$, and $\{G_n\}$ be a sequence of $(d,\varepsilon)$-locally sparse $d$-regular graphs on $[n]$.
Then we have ${\rm{th}}(G)=n^{-\frac{2}{d}}$.
\end{corollary}

Now, we show that Corollary~\ref{cor:threshold_of_locally_sparse_regular_graph} is almost tight.
Let $d,n\in\mathbb{N}$ and $0<\varepsilon<1/2$, where $d$ divides $n$.
We give a construction for $n$-vertex $d$-regular graph $G$ such that the minimum edge boundary over all subsets $U\subseteq[n]$ of size $d\le |U|\le{\varepsilon}n$ is $d$ and has ${\rm{th}}(G)\gg{n^{-2/d}}$.
This implies the lower bound of edge boundary in the definition of $(d,\varepsilon)$-locally sparse graphs is the best possible for $d$-regular graphs.
Let $K_i$ be an $d$-clique on vertex set $\left\{(i-1)d+1,\ldots,id\right\}$, for every $i=1,\ldots,n/d$.
When $d$ is an even number, $G$ is obtained from the disjoint union $\bigcup^{n/d}_{i=1}K_i$ by adding edges $\left\{(i-1)d+\frac{d}{2}+j,id+j\right\}$ for $i=1,\ldots,n/d,\text{ }j=1,\ldots,d/2$, where $n+j=j$.
When $d$ is an odd number, $G$ is obtained from the disjoint union $\bigcup^{n/d}_{i=1}K_i$ by adding edges $\left\{(i-1)d+\frac{d-1}{2}+j,id+j\right\}$ and $\left\{(i-1)d+d,(\frac{n}{2d}+i-1)d+d\right\}$ for $i=1,\ldots,n/d,\text{ }j=1,\ldots,(d-1)/2$, where $n+j=j$.
It is not hard to check that the resulting graph $G$ is $d$-regular and the minimum edge boundary of $G$ over all subsets $U\subseteq[n]$ of size $d\le |U|\le{\varepsilon}n$ is $d$. Since $G$ contains a $K_d$-factor, we have ${\rm{th}}(G)\ge{\left(n^{-1}\log^{1/(d-1)} n\right)^{2/d}}$.
On the other hand, Proposition~\ref{prop:Kelly_edge_spread} implies ${\rm{th}}(G)\le{n^{-2/d}\log n}$.

\medskip
\noindent\textbf{Remarks.}
\begin{itemize}
\item The upper bound $n^{-1/d}$ of ${\rm th}(G)$ is best possible if $G$ is asymptotically maximal subject to the $(d,\alpha)$-degenerate. Indeed, if $G$ satisfies $e(G)\ge dn - O(1)$, then we have $p_E(G)\ge n^{-1/d}$ by computing the expected number of copies of $G$.

\medskip
\item {Our result on thresholds improves on Proposition~\ref{prop:Kelly_edge_spread} for $(d,\alpha)$-degenerate graphs by eliminating the logarithmic factor.
This contributes to a better understanding of the ratio ${\rm th}(G)/p_E(G)$ (see Conjecture \ref{conj:conj_seocnd_kk}) -- all known results on this ratio is at most $\log n$ and Theorem \ref{thm:asthm0} shows that this gap vanishes for $(d,\alpha)$-degenerate $k$-graphs $G$ with $e(G)\ge dv(G) - O(1)$.
This question is also studied in a recent paper by Balogh, Bernshteyn, Delcourt, Ferber, and Pham in~\cite{Balogh2024SunflowersIS} for hypergraphs of bounded VC-dimension.
}

\medskip
\item
Riordan's result~\cite{Riordan2000SpanningSO} implies that ${\rm{th}}(G_n)\le n^{-1/d}$ for $(d, d)$-degenerate graphs.
Thus, our result improves on Riordan's result by allowing arbitrarily small $\alpha>0$.
\medskip
\item In order to get rid of the logarithmic term from the threshold, the condition $\alpha>0$ in the definition of $(d,\alpha)$-degenerate hypergraphs is best possible in the following sense: if $G_n$ is a Hamilton cycle or a spanning tree, then it is $(1,0)$-degenerate and has a threshold $\log n/n$ (see also the examples above).

\end{itemize}

\section{Spreadness} \label{sec::spr}

Recently, random graph theory was revolutionized by Frankston, Kahn, Narayanan, and Park's~\cite{FKNP2021} proof of the \emph{fractional expectation threshold vs.~threshold} conjecture of Talagrand~\cite{Talagrand2010AreMS} and Park and Pham's~\cite{Park2023APO} proof of the even stronger (abstract version of) Kahn--Kalai conjecture~\cite{kahn2007thresholds}.
A crucial corollary of Talagrand's conjecture relates the so-called \emph{spread measure} with thresholds.
Roughly speaking, for the study of random graphs, this connection implies that if one can find a probability measure on the copies of $G$ in $K_n$ with good \emph{spread}, then $G$ w.h.p. appears in $G(n,p)$.

\begin{definition}\label{def:edge_spread}
Let $q\in[0,1]$. Assume that $\mathcal{H}$ is a hypergraph on the vertex set $V$,  and $\mu$ is a probability measure on the edge set of $\mathcal{H}$.
We say $\mu$ is \emph{$q$-spread} if for every $S\subseteq{V}$, we have ~\[\mu\left(\{A\in{E(\mathcal{H})}:S\subseteq{A}\}\right)\le{q^{|S|}}.\]

\end{definition}

Pham, Sah, Sawhney, and Simkin introduced a notation of \emph{vertex spreadness} in~\cite{pham2022toolkit}, where they consider a random embedding of $T\rightarrow{H}$ and analyze the probability that a subsequence of vertices of $T$ is mapped to a subsequence of vertices of $H$.
Kelly, M{\"u}yesser, and Pokrovskiy~\cite{kelly2023optimal} extended their result to the general settings.
\begin{definition}[Vertex-spread]
Let $X$ and $Y$ be two finite sets, and let $\mu$ be a probability distribution over injections $\psi: X\rightarrow{Y}$. For $q\in[0,1]$, we say that $\mu$ is $q$-\emph{vertex spread} if for every $s \le |X|$ and every pair of subsequences of distinct vertices $x_1,\dots,x_s\in{X}$ and $y_1,\dots,y_s\in{Y}$,
\[\mu\left(\{\psi(x_i)=y_i,\forall{i\in[s]}\}\right)\le{q^s}.\]

\end{definition}

A \emph{hypergraph embedding}~$\psi: G\rightarrow{H}$ of a hypergraph~$G$ into a hypergraph~$H$~is an injective map~$\psi: V(G)\rightarrow{V(H)}$ which maps edges of $G$ to edges of $H$, so there is an embedding of $G$ into $H$ if and only if $H$ contains a subgraph isomorphic to $G$.
\begin{fact}\label{fact:trivial_vertex_spread}
Let $G$ be an $n$-vertex $k$-uniform hypergraph. Then the uniform distribution $\mu$ on embeddings $\psi$ of $G$ in $K^{(k)}_n$ is $(e/n)$-vertex spread.
\end{fact}

\begin{proof}
Fixing $s\ge 0$ and vertices $x_1,\dots, x_s\in V(G)$ and $y_1,\dots, y_s\in V(K_n^{(k)})$, we have
\[
\mu\left(\{\psi: \psi(x_i)=y_i,\,\forall{i\in[s]}\}\right) = \frac{(n-s)!}{n!}\le (e/n)^s. \qedhere
\]
\end{proof}

The following result is proved by Kelly, M{\"u}yesser, and Pokrovskiy in~\cite{kelly2023optimal} which allows us to connect spread distributions and vertex spread distributions.
{Moreover, combined with the main theorem in~\cite{FKNP2021}, it gives an upper bound for ${\rm{th}}(G)$.}

\begin{proposition}[Proposition 1.17 in \rm{\cite{kelly2023optimal}}]\label{prop:Kelly_edge_spread}
Let $k,\Delta\in{\mathbb{N}}$ and $C$ be a positive constant, there exists a constant $C_{\ref{prop:Kelly_edge_spread}}=C_{\ref{prop:Kelly_edge_spread}}(C,k,\Delta)>0$ such that the following holds for all sufficiently large $n$.
Let $H$ and $G$ be two $n$-vertex $k$-uniform hypergraphs.
Let $\mathcal{H}$ be the set of subgraphs of $H$ which are isomorphic to $G$ and suppose there is a $(C/n)$-vertex spread distribution $\mu$ on embeddings $\psi:G\rightarrow{H}$ and $\Delta(G)\le{\Delta}$.
Then there exists a $\left(C_{\ref{prop:Kelly_edge_spread}}/n^{1/m_1(G)}\right)$-spread distribution $\mu^{\prime}$ on $\mathcal{H}$ such that for any subset $S\subseteq{E(H)}$ with $v$ vertices and $c$ components, we have $\mu^{\prime}\left(\{A\in{\mathcal{H}}:S\subseteq{A}\}\right)\le{(k\Delta C^2/n)^{v-c}}$.
In particular, we have ${\rm{th}}(G)\le{n^{-1/m_1(G)}\log n}$.
\end{proposition}

{We first state the following stronger notion of spreadness introduced by Spiro~\cite{spiro2023smoother}.
We remark that a similar notion of ``superspread'' was also introduced by Espuny D\'{\i}az and Person~\cite{Daz2023SpanningI}.}

\begin{definition}\label{def:smoother-spread}
Let $q\in[0,1]$, and let~$r_0,r_1,\dots,r_l\in\mathbb{N}$ be a decreasing sequence of positive integers.
Let $\mathcal{H}$ be an $r_0$-bounded hypergraph on $V$, and let $\mu$ be a probability measure on $E(\mathcal{H})$. We say $\mu$ is $(q;r_0,\dots,r_l)$\emph{-spread} if the following holds for all $i\in[l]$:
\[\mu\big(\{A\in{E(\mathcal{H})}~:~|A\cap{S}|\ge t\}\big)\le{q^t}~\text{for~all~}t\in[r_i,r_{i-1}]\text{~and~}S\in{\bigcup_{j=r_i}^{r_{i-1}}\Bigg\{\binom{A}{j}~:~A\in{E(\mathcal{H})}\Bigg\}}.\]

\end{definition}

The following result is proved by Spiro~\cite{spiro2023smoother} and is also a stronger version of the main theorem in~\cite[Theorem 1.6]{FKNP2021}.

\begin{theorem}[\rm{Theorem 1.3 in \cite{spiro2023smoother}}]\label{thm: stronger-fknp}
There exists an absolute constant $K_{\ref{thm: stronger-fknp}}>0$ such that the following holds for all  $C\ge{K_{\ref{thm: stronger-fknp}}}$.
Let $r_0,r_1,\dots,r_l\in\mathbb{N}$ be a decreasing sequence of positive integers, where $r_l=1$, and $\mathcal{H}$ be an $r_0$-bounded hypergraph on $V$. If there is a $(q;r_0,\dots,r_l)$-spread distribution on $E(\mathcal{H})$, then a set $W$ of size $C\ell q|V|$ chosen uniformly at random from $V$ contains an edge of $\mathcal{H}$ with probability at least $1-K_{\ref{thm: stronger-fknp}}/(Cl)$.
\end{theorem}

Now we are ready to state our main result, improving Proposition \ref{prop:Kelly_edge_spread}.
Indeed, we show that for $(d,\alpha)$-degenerate $k$-uniform hypergraph $G$, the distribution $\mu'$ given by Proposition~\ref{prop:Kelly_edge_spread} satisfies Spiro spreadness with constant $l$, and thus obtain a saving of the logarithmic factor.

\begin{theorem}\label{thm:asthm}
Let $k,\Delta\in\mathbb{N}$,~$d,\alpha,C^{\prime}>0$, there exists a constant $C>0$ such that the following holds for all sufficiently large $n$.
Let $G$ be an $n$-vertex $(d,\alpha)$-degenerate $k$-uniform hypergraph with $\Delta(G)\le\Delta$ and let $H$ be an $n$-vertex $k$-uniform hypergraph.
If there is a $(C^{\prime}/n)$-vertex spread distribution $\mu$ on embeddings $\psi:G\rightarrow{H}$,
then there exists a $\left(C/n^{1/d};r_0,\ldots,r_l\right)$-spread distribution $\mu^{\prime}$ on copies of $G$ in $H$
where $r_0=dn,r_1=\varepsilon n/k,r_i=r_{i-1}/{n^{9\alpha/(10d)}}$ for $2\le i\le{l-1}$, and $r_l=1$.
In particular, we have ${\rm{th}}(G)\le{n^{-1/d}}$.
\end{theorem}

Clearly, Theorem \ref{thm:asthm} implies Theorem \ref{thm:asthm0}.
Moreover, we shall include a slight refinement of Theorem \ref{thm:asthm} in Section 4, using a recent work of Kelly~\cite{kelly2024spread}.

With additional (and separate) work on the vertex spreadness, Theorem~\ref{thm:asthm} can be used to obtain \emph{robust thresholds} for $(d,\alpha)$-degenerate hypergraphs, similar to some results in~\cite{kelly2023optimal} which was derived from establishing vertex spreadness in host graphs and applying Proposition~\ref{prop:Kelly_edge_spread}.
One prominent such example is the robust threshold of $r$-th power of Hamilton cycles in graphs, which was recently determined by Joos--Lang--Sanhueza-Matamala~\cite{joos2023robusthamiltonicity}.

\section{Proof of Theorem~\ref{thm:asthm}}\label{sec:section_3}

We first include two auxiliary results.
The first one is proved in~\cite{kelly2023optimal}.

\begin{lemma}\label{lem:isomorphism-counting}{\rm(Lemma 6.1~in~\cite{kelly2023optimal})}
Let $G$ and $H$ be $n$-vertex $k$-uniform hypergraphs.
If $F\subseteq{H}$ has $v$ vertices and $c$ components, then there are at most
\[n^c\big(k\Delta(G)\big)^{v-c}\]
hypergraph embeddings $G[X]\rightarrow{H[V(F)]}$ where $X\in{\binom{V(G)}{v}}$ and $F\subseteq{H[\psi(X)]}$.
\end{lemma}

The second one is a simple counting result and a specialized version was proved in~\cite{kelly2023optimal}.
\begin{lemma}\label{lem:subgraph-counting}
Let $G$ be an $n$-vertex $k$-uniform hypergraph with $\Delta(G)\le\Delta$.
If $S$ is a nonempty subset of $E(G)$, then the number of subgraphs of $G$ with $t$ edges, no isolated vertices, all in $S$, and $c$ components is at most
\begin{equation}\label{arabic:subgraph_counting}
2^{(k+1)t}{{\Delta}^{t}}\binom{k|S|}{c}.
\end{equation}
\end{lemma}

\begin{proof}
Let $F$ be a subgraph of $G$ which is induced by the $t$ edges in $S$. Let $F_1,F_2,\dots,F_c$ be the $c$ components of $F$, and $v_1,v_2,\dots,v_c$ be a rooted vertex of $F_1,F_2,\dots,F_c$, respectively.
We have
\begin{enumerate}
        \item[$\bullet$] the number of choices of $v_1,v_2,\dots,v_c$ is at most $\binom{k|S|}{c}$;
        \item[$\bullet$] for each $i\in[c]$, let $|E(F_i)|=t_i$, then the number of choices of the equation $t=\sum_{i=1}^ct_i$ is at most $\binom{t-1}{c-1}$;
        \item[$\bullet$] for each $i\in[c]$, the number of choices of $F_i\subseteq{S}$ is at most $2^{kt_i}{\Delta}^{t_i}$, since there are at most $2^{kt_i}$ degree sequences $(d_1, \dots, d_{n_i})$ for $F_i$, where $kt_i=\sum_{j=1}^{n_i}d_j$.
\end{enumerate}
Thus, we have the number of choices of $F$ is at most $2^{(k+1)t}{{\Delta}^{t}}\binom{k|S|}{c}$.
\end{proof}

Now we are ready to prove Theorem~\ref{thm:asthm}.
As mentioned in Section~\ref{sec::spr}, we are going to show that for $(d,\alpha)$-degenerate $k$-uniform hypergraph $G$, the distribution $\mu'$ given by Proposition~\ref{prop:Kelly_edge_spread} satisfies Spiro spreadness with some constant $l$.

\begin{proof}[Proof of Theorem~\ref{thm:asthm}]
Let $1/n\ll{1/\tilde{C}}\ll{1/C}\ll{1/\Delta}<1$, and $d,\alpha>0$.
Let $H$ be an $n$-vertex $k$-uniform hypergraph and $G$ be an $n$-vertex $(d,\alpha)$-degenerate $k$-uniform hypergraph with $\Delta(G)\le\Delta$.
Assume $\mu$ is a $(C/n)$-vertex spread distribution on embeddings $\psi$ of $G$ in $H$.
Let $\mathcal{H}$ be the $(dn)$-bounded hypergraph with vertex set $E(H)$ and edge set the set of subgraphs of $H$ which are isomorphic to $G$.
By Proposition~\ref{prop:Kelly_edge_spread}, the distribution $\mu^{\prime}$ on $E(\mathcal{H})$ is $\left(C_{\ref{prop:Kelly_edge_spread}}/n^{1/m_1(G)}\right)$-spread.
Since $G$ is $(d,\alpha)$-degenerate, there exists a constant $0<\varepsilon<1$ such that for every connected subgraph of $G$ with $v$ vertices and $t$ edges, $0<t\le{\varepsilon}n/k$, we have $v\ge\frac{t}{d}+1+\frac{\alpha}{d}$.

Now we show that $\mu^{\prime}$ is $\left(\frac{\tilde{C}}{n^{1/d}};dn,\varepsilon n/k,\varepsilon n^{1-9\alpha/(10d)}/k,\dots,\varepsilon n^{9\alpha/(10d)}/k,1\right)$-spread.
To that end, we assume $t\in[dn]$ and $S\subseteq{E(H)}$, $S\neq\emptyset$.
For an edge set $S\subseteq{E(H)}$, we let $c(S)$ be the number of components that are induced by $S$ and $v(S)$ be the number of vertices contained in $S$.
Let $F$ be the subgraph of $H$ with edge set $S$ and subject to that, the fewest number of vertices.
We can assume that $F$ is isomorphic to a subgraph of $G$, since otherwise $\mu^{\prime}\left(\{A\in{E(\mathcal{H})}:F\subseteq{A}\}\right)=0$.

If $t\in[\varepsilon n/k,dn]$, since $m_1(G)\le{d}$, and $\mu^{\prime}$~is $(C_{\ref{prop:Kelly_edge_spread}}n^{-1/m_1(G)})$-spread, then we have
\begin{align*}
\mu^{\prime}\big(\{A\in{E(\mathcal{H})~:~|A\cap{S}|\ge t}\}\big)\le {} &\sum_{T\in\binom{S}{t}}\mu^{\prime}\big(\{A\in{E(\mathcal{H})~:~T\subseteq A}\}\big) \\
\le {} & \binom{|S|}{t}\left(\frac{C_{\ref{prop:Kelly_edge_spread}}}{n^{1/m_1(G)}}\right)^t \\
\le {} & {\left(\frac{e|S|}{t}\right)^t}\left(\frac{C_{\ref{prop:Kelly_edge_spread}}}{n^{1/d}}\right)^t\le\left(\frac{(ekd/\varepsilon)C_{\ref{prop:Kelly_edge_spread}}}{n^{1/d}}\right)^t,
\end{align*}
as desired.

Next, we consider the case $|S|, t\le \varepsilon n/k$.
For any $c\le t$, we claim that
\begin{equation}\label{arabic:|S|_choose_c}
        \binom{k|S|}{c}\le\left\{
        \begin{array}{rl}
        {n^{(9\alpha c)/(10d)}} & \text{if }|S|,t\in[1,\varepsilon n^{9\alpha/(10d)}/k], \\
        {2^tn^{\alpha c/d}}      & \text{otherwise}.
        \end{array}\right.
\end{equation}
Indeed, when $|S|,t\in[1,\varepsilon n^{9\alpha/(10d)}/k]$, we have $\binom{k|S|}{c}\le (k|S|/c)^c \le{(n^{9\alpha/(10d)})^c}$.
Moreover, for the case of $|S|,t\in[\varepsilon n^{1-9i\alpha/(10d)}/k,\varepsilon n^{1-9(i-1)\alpha/(10d)}/k]$, where $i\in[l]$ and $l=\lceil10d/(9\alpha)\rceil$,
if $c\ge{\frac{ek|S|}{n^{\alpha/d}}}$, we have $\binom{k|S|}{c}\le{(n^{\alpha/d})^c}$; otherwise, we have $c<{\frac{ek|S|}{n^{\alpha/d}}}\le{t/n^{\alpha/(20d)}}$, which implies $\binom{k|S|}{c}\le{(k|S|)^c}\le{n^{(t/\log n)}}\le{2^t}$, since $k|S|<{\varepsilon n}<n$.

Recall that $G$ is $(d,\alpha)$-degenerate, then for every connected subgraph of $G$ with $v$ vertices and $t$ edges, $0<t\le{\varepsilon}n/k$, we have $v\ge\frac{t}{d}+1+\frac{\alpha}{d}$.
If $T$ is a subgraph of $G$ with $c$ components $T_1,\ldots,T_c$ such that each component of $T$ has size at most $\varepsilon n/k$, then we have $v(T)=\sum^c_{i=1}v(T_i)\ge\sum^c_{i=1}\left(e(T_i)/d+1+\alpha/d\right)={\frac{t}{d}+(1+\frac{\alpha}{d})c}$.

By Lemma~\ref{lem:isomorphism-counting}, Lemma~\ref{lem:subgraph-counting}, and Proposition~\ref{prop:Kelly_edge_spread}, we have
\begin{align*}
\mu^{\prime}\left(\{A\in{E{(\mathcal{H})}}:|A\cap{S}|\ge t\}\right)\le {}&\sum_{X\in\binom{S}{t}}\mu^{\prime}\left(\{A\in{E{(\mathcal{H})}}:X\subseteq{A}\}\right) \\
{\le} {}&\sum_{X\in\binom{S}{t}}\left(\frac{k\Delta C^2}{n}\right)^{v(X)-c(X)} \\
\le {}&\sum^{t}_{c=1}\sum_{X\in\binom{S}{t}:c(X)=c}\left(\frac{k\Delta C^2}{n}\right)^{v(X)-c} \\
\overset{(\ref{arabic:subgraph_counting})}{\le} {} &\sum_{c=1}^t2^{(k+1)t}{{\Delta}^{t}}\binom{k|S|}{c}\left(\frac{k\Delta C^2}{n}\right)^{\frac{t}{d}+\frac{\alpha c}{d}} \\
\overset{(\ref{arabic:|S|_choose_c})}{\le} {} &\sum_{c=1}^t2^{(k+1)t}{\Delta}^t 2^tn^{\alpha c/d}\left(\frac{k\Delta C^2}{n}\right)^{\frac{t}{d}+\frac{\alpha c}{d}} \\
\le {}&t\cdot \left(\Delta 2^{k+2}\right)^t\left(\frac{k\Delta C^2}{n}\right)^{\frac{t}{d}} \\
\le {}&{\left(\frac{\left(2^{k+3}(kC^2\Delta^{d+1})^{1/d}\right)}{n^{1/d}}\right)^t},
\end{align*}
where the second inequality follows from Proposition~\ref{prop:Kelly_edge_spread}.

In particular, when $|S|,t\in[1,\varepsilon n^{9\alpha/(10d)}/k]$, we obtain that
\begin{align}
\mu^{\prime}\left(\{A\in{E{(\mathcal{H})}}:|A\cap{S}|\ge t\}\right)\le {} &\sum_{c=1}^t2^{(k+1)t}{{\Delta}^{t}}\binom{k|S|}{c}\left(\frac{k\Delta C^2}{n}\right)^{\frac{t}{d}+\frac{\alpha c}{d}} \nonumber \\
\overset{(\ref{arabic:|S|_choose_c})}{\le} {} &\sum_{c=1}^t2^{(k+1)t}{\Delta}^t n^{0.9\alpha c/d}\left(\frac{k\Delta C^2}{n}\right)^{\frac{t}{d}+\frac{\alpha c}{d}} \nonumber \\
\le {}& t\cdot n^{-0.1\alpha/d}\left(\Delta 2^{k+1}\right)^t\left(\frac{k\Delta C^2}{n}\right)^{\frac{t}{d}} \nonumber \\
\le {}&n^{-0.1\alpha/d}{\left(\frac{\left(2^{k+2}(kC^2\Delta^{d+1})^{1/d}\right)}{n^{1/d}}\right)^t}. \label{eq:111}
\end{align}
Therefore, we obtain that $\mu^{\prime}$ is $\left(\tilde{C}/n^{1/d};dn,\varepsilon n/k,\varepsilon n^{1-9\alpha/(10d)}/k,\dots,\varepsilon n^{9\alpha/(10d)}/k,1\right)$-spread,
where $\tilde{C}=2^{(k+3)}ekd{C_{\ref{prop:Kelly_edge_spread}}}{(kC^2\Delta^{d+1})^{\frac{1}{d}}}/\varepsilon$.

Note that $l=\max\{2,\lceil{(10d)/(9\alpha)}\rceil\}$, we let $C_{\ref{thm:asthm}}\ge{(K_{\ref{thm: stronger-fknp}}\tilde{C}l)}$ be a positive constant.
When $p\ge{C_{\ref{thm:asthm}}/n^{1/d}}$, by Theorem~\ref{thm: stronger-fknp}, $H_p$ contains a copy of $G$ with probability at least $1-o_{C_{\ref{thm:asthm}}}(1)$, where $H_p$ denotes the random subset obtained by keeping every edge of $H$ independently with probability $p$.
\end{proof}

\section{Additional note}

When we were finalizing this project, we heard from Tom Kelly that
he introduced yet another strengthening of the spreadness~\cite{kelly2024spread}, which can boost our result to a ``semi-sharp'' threshold.
\begin{definition}\label{def:sharper_spread}
Let $q,\gamma\in[0,1]$, and let $r,r^{\prime}\in\mathbb{N}$ be positive integers such that $r>r^{\prime}$, and let $\mathcal{H}$ be a hypergraph on vertex set $V$.
We say the probability measure $\mu$ on $E(\mathcal{H})$ is $(q,\gamma;r,r^{\prime})$-\emph{spread} if the following holds:
\[\mu\big(\{A\in{E(\mathcal{H})}~:~|A\cap{S}|\ge t\}\big)\le{\gamma q^t}~\text{for~all~}t\in[r^{\prime},r]\text{~and~}S\in{\bigcup_{j=r^{\prime}}^r\Bigg\{\binom{A}{j}~:~A\in{E(\mathcal{H})}\Bigg\}}.\]
\end{definition}

The following is a slightly stronger version of the result of Spiro~\cite[Theorem 2.7]{spiro2023smoother} that incorporates Definition~\ref{def:sharper_spread}.
\begin{theorem}[\cite{kelly2024spread}]\label{thm:adapted_Spiro_spread_fknp}
Let $q,\gamma\in[0,1]$, and $r_0,r_1,\dots,r_{\ell}\in\mathbb{N}$ be a decreasing sequence of positive integers, where $r_{\ell}=1$, and $\mathcal{H}$ be an $r_0$-bounded hypergraph on $V$, and let $C_{\ref{thm:adapted_Spiro_spread_fknp}}\ge 8$.
If there is a $(q;r_0,\dots,r_{\ell})$-spread distribution on $E(\mathcal{H})$ which is $(q,\gamma;r_{\ell-1},1)$-spread and $W$ is a set of size $2{\ell}C_{\ref{thm:adapted_Spiro_spread_fknp}}q|V|$ chosen uniformly at random from $V$, then $W$ contains an edge of $\mathcal{H}$ with probability at least
\[1-6\ell^2\left(\frac{C_{\ref{thm:adapted_Spiro_spread_fknp}}}{4}\right)^{-r_{\ell-1}/2}-2\exp\left(-\frac{\ell C_{\ref{thm:adapted_Spiro_spread_fknp}}q|V|}{4}\right)-\frac{8\gamma}{C_{\ref{thm:adapted_Spiro_spread_fknp}}\ell}.\]
\end{theorem}
For the negative terms in the probability above, the first term is exponential in $r_{\ell-1}$, the second one is exponential in $|V|$, and the last is linear in $\gamma/\ell$.
Therefore, if we have $r_{\ell-1}\to \infty$, $|V|\to \infty$ and $\gamma/\ell \to 0$, then the claimed probability above is $1-o(1)$.

Using Theorem~\ref{thm:adapted_Spiro_spread_fknp}, we can improve the error probability in our result and obtain a ``semi-sharp'' threshold for $(d,\alpha)$-degenerate hypergraphs.
Note that this also implies that \emph{all nearly maximally $(d,\alpha)$-degenerate hypergraphs have sharp thresholds}.

\begin{theorem}\label{thm:asthm++}
Let $k,\Delta\in\mathbb{N}$,~$d,\alpha,C^{\prime}>0$, there exists a constant $C, C^*>0$ such that the following holds for all sufficiently large $n$.
Let $G$ be an $n$-vertex $(d,\alpha)$-degenerate $k$-uniform hypergraph with $\Delta(G)\le\Delta$ and let $H$ be an $n$-vertex $k$-uniform hypergraph.
If there is a $(C^{\prime}/n)$-vertex spread distribution $\mu$ on embeddings $\psi:G\rightarrow{H}$,
then there exists a $\left(C/n^{1/d};r_0,\ldots,r_l\right)$-spread distribution $\mu^{\prime}$ on copies of $G$ in $H$ which is $(C/n^{1/d},n^{-\alpha/(10d)};r_{\ell-1},1)$-spread,
where $r_0=dn,r_1=\varepsilon n/k,r_i=r_{i-1}/{n^{9\alpha/(10d)}}$ for $2\le i\le{l-1}$, and $r_l=1$.
Moreover, if $p\ge C^*n^{-1/d}$, then $P[G\subseteq{G^{(k)}(n,p)}] = 1-o_n(1)$.
\end{theorem}

Note that our proof of Theorem~\ref{thm:asthm} gives that $\mu'$ is $(C/n^{1/d},n^{-\alpha/(10d)};r_{\ell-1},1)$-spread by~\eqref{eq:111}.
Then Theorem~\ref{thm:asthm++} follows from Theorem~\ref{thm:adapted_Spiro_spread_fknp}.

\section*{Acknowledgment}
This work was partially done when the third author visited the Beijing Institute of Technology in October 2023. Jie Han is partially supported by National Natural Science Foundation of China (12371341). Haoran Luo is partially supported by UIUC Campus Research Board 24012 and NSF grant FRG-DMS 2152488.

We thank Yury Person, Olaf Parczyk for valuable discussions, and Tom Kelly for suggesting Theorem 4.3 and sharing his manuscript~\cite{kelly2024spread} with us.

\end{document}